\documentclass[11pt,reqno]{amsart}
\usepackage{amssymb,latexsym}
\usepackage{graphicx}
\usepackage{comment}
\usepackage{url}		%does nice formatting of URLs
\usepackage{color}
\usepackage{hyperref}
\usepackage{url}
\usepackage{breakurl}
\newcommand{\bburl}[1]{\textcolor{blue}{\url{#1}}}

\calclayout

\makeatletter
\newcommand{\monthyear}[1]{%
  \def\@monthyear{\uppercase{#1}}}
\newcommand{\volnumber}[1]{%
  \def\@volnumber{\uppercase{#1}}}
\AtBeginDocument{%
\def\ps@plain{\ps@empty
  \def\@oddfoot{\@monthyear \hfil \thepage}%
  \def\@evenfoot{\thepage \hfil \@volnumber}}
\def\ps@firstpage{\ps@plain}
\def\ps@headings{\ps@empty
  \def\@evenhead{%
    \setTrue{runhead}%
    \def\thanks{\protect\thanks@warning}%
    \uppercase{The Fibonacci Quarterly}\hfil}%
  \def\@oddhead{%
    \setTrue{runhead}%
    \def\thanks{\protect\thanks@warning}%
    \hfill\uppercase{Hypergeometric Template}}%
  \let\@mkboth\markboth
  \def\@evenfoot{%
    \thepage \hfil \@volnumber}%
  \def\@oddfoot{%
    \@monthyear \hfil \thepage}%
  }%
\footskip=25pt
\pagestyle{headings}%
}
\makeatother

\newcommand{\Q}{{\mathbb Q}}

\theoremstyle{plain}
\numberwithin{equation}{section}

\begin{document}
%% replace the values in the next three lines by the correct information
\monthyear{Month Year}
\volnumber{Volume, Number}
\setcounter{page}{1}

\title{Geometric Proofs of the Irrationality of Square-Roots for Select Integers}
\author{Zongyun Chen, Steven J. Miller, Chenghan Wu}

\address{193/198, Soi 9, Koolpunt Ville 7, Moo 4, Muang Chiangmai District, Chiangmai Province, Thailand}
\email{marinachenzongyun@gmail.com}

\address{Department of Mathematics and Statistics\\
                Williams College\\
                Williamstown, MA\\}
\email{sjm1@williams.edu}

\address{9 Jixing Road, Yubei District, Chongqing, China}
\email{2907360@qq.com}

\thanks{This paper is an outgrowth of a talk given by the second named author and attended by the other authors at the Math League International Summer Tournament in Summer 2024 in Trenton, NJ. It is a pleasure to thank the organizers and participants for creating such a lively atmosphere, especially Dan Flegler, John Hagen, Rui Hu, and Adam Raichel.}

 \maketitle
%%%%%%%%%%%%%%%%%%%%%%%%%%%%%%%%%%%%%%%%%%%%%%%%%%%%%%%%%%%%%%%%%%%%%%%%%%%%%%%%%%%%%%
%%%%%%%%%%%%%%%%%%%%%%%%%%%%%%%%%%%%%%%%%%%%%%%%%%%%%%%%%%%%%%%%%%%%%%%%%%%%%%%%%%%%%%
%%%%%%%%%%%%%%%%%%%%%%%%%%%%%%%%%%%%%%%%%%%%%%%%%%%%%%%%%%%%%%%%%%%%%%%%%%%%%%%%%%%%%%
\section{Introduction}

The positive integers $1, 2, 3, \dots$ are not surprisingly one of the most important sequences in mathematics, and typically the first encountered. Quickly one meets interesting sub-sequences, such as the primes ($2, 3, 5, 7, 11, \dots$), the perfect squares ($1, 4, 9, 16, 25, \dots$) and the Fibonacci numbers ($1, 2, 3, 5, 8, \dots$) to name just a few. These are well studied and arise in numerous places; see the On-line Encyclopedia of Integer Sequences \cite{OEIS} for details and properties of these and others.

Almost all integers have irrational square-roots, with the percent of $n \le x$ with $\sqrt{n} \not\in\Q$ approximately $100 \cdot x^{-1/2}\%$. The standard proof uses the property that if a prime $p$ divides a product $xy$ then $p|x$ or $p|y$ or both (see for example \cite{MS} for a proof) and the Fundamental Theorem of Arithmetic (every integer can be written uniquely as a product of primes in increasing order; see \cite{HW}). \begin{quote} Assume a non-square $n > 1$ has a rational square-root; thus we can write $\sqrt{n} = a/b \in \Q$ with $a, b$ relatively prime integers and without loss of generality it suffices to consider $n$ that are square-free, as if $n = m_1 m_2^2$ then $\sqrt{n} = \sqrt{m_1} \cdot m_2$. Then $nb^2 = a^2$. As $n > 1$ is square-free, there is a prime $p$ that divides $n$. Thus $p|a^2$ so $p|a$ and we can write $a$ as $\alpha p$. Substituting yields $nb^2 = \alpha^2 p^2$; as $n$ is square-free and a multiple of $p$, we must have $n/p$ is an integer relatively prime to $p$ and thus $p|b^2$. A similar argument now shows $b = \beta p$, contradicting $a$ and $b$ are relatively prime and therefore $\sqrt{n}$ is irrational.  \end{quote}

There's a lot of interesting history on this proof; if we don't use the property that if a prime divides a product then it divides at least one factor, we can mimic the above argument, but only by essentially reproving the result case by case. For example, if $n=2$ then we would have $2b^2 = a^2$. If $a = 2\alpha+1$ is odd then $a^2 = 4\alpha^2 + 4\alpha + 1$ is odd, and thus cannot be a multiple of 2, and thus $a =2 \alpha$. Similarly if $n=3$ we would have $3b^2 = a^2$ and 3 must divide the right hand side as it divides the left. We can write $a = 3\alpha+r$ with $r \in \{0, 1, 2\}$ and note $$a^2\ =\ 9 \alpha^2 + 6\alpha r + r^2 \ = \ 3(3\alpha^2 + 2 \alpha r) + r^2,$$ which is only a multiple of 3 when $r=0$ as $1^2=1$ and $2^2 = 4$ are not multiples of 3. Do we need the fact that if a prime divides a product it divides at least one factor? It turns out we can verify enough of this property for any specific prime, in particular we can show if $p$ divides a number of the form $mp+r$ for $0 \le r < p$ then $r$ must be zero; doing so just requires showing $r^2$ is not a multiple of $p$ for $1 \le r \le p-1$. Unfortunately we seem to need the product property to argue generally.
%for any $n$ that is square-free, we can always directly compute $r^2$ and show that it is not a multiple of $n$.

We thus know that the elements of $\mathcal{I} := \{2, 3, 5, 6, 7, 8, 10, 11, 12, 13, 14, 15, 17, \dots\}$ have irrational square-roots, provable using the above methods. What if we wish to avoid using a prime dividing a product divides at least one factor? Which elements in $\mathcal{I}$ can be proved to have irrational square-roots using geometric methods? To be fair Euclid proved in Book VII of his classic Elements the result that if $p|xy$ then $p|x$ or $p|y$, using of course different language than we use today;\footnote{From the Wikipedia entry on the Fundamental theorem of arithmetic \cite{Wi}: \emph{If two numbers by multiplying one another make some number, and any prime number measure the product, it will also measure one of the original numbers.}} thus by geometric argument we mean using areas of figures.

A classical geometric proof for the irrationality of $\sqrt{2}$ was introduced by Stanley Tennenbaum \cite{Tennenbaum} in the 1950s and later popularized by John H. Conway in his article \textit{The Power of Mathematics} \cite{Conway}. For completeness we quickly review their proof, paraphrasing from \cite{MillerMontague1}. \begin{quote}
Suppose that $(a/b)^2 = 2$ for integers $a$ and $b$; we may assume $a$ and $b$ are the smallest such numbers. Thus $a^2 = 2b^2$, which we interpret geometrically as the area of two squares of side length $b$ that equals the area of one square of side length $a$. Thus, if we consider Figure \ref{fig:sqrt2}, \begin{figure}[h]
\begin{center}
\scalebox{.7}{\includegraphics{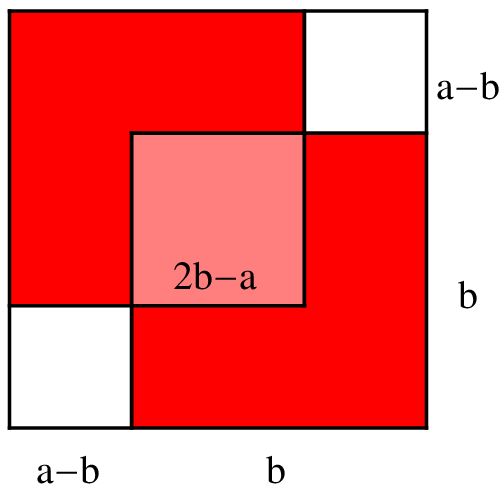}}
\caption{\label{fig:sqrt2} Geometric proof of the irrationality of $\sqrt{2}$.}
\end{center}\end{figure}  the total area covered by the squares with side length $b$ (double counting the overlapping, pink region) equals the area of the larger square with side length $a$. Therefore the pink, doubly counted part, which is a square of side length $2b-a$, has area equal to that of the two white, uncovered squares (each of side length $a-b$). Therefore $(2b-a)^2 = 2(a-b)^2$ or $\sqrt{2} = (2b-a)/(a-b)$. We now show this implies we have $\sqrt{2}$ equal to the ratio of two smaller integers, a contradiction completing the proof. Clearly $2b - a > 0$, as if not we would have $2b \le a$ and thus $4b^2 \le a^2 = 2b^2$. For the upper bound, note $2b - a < a$, as otherwise $2b-a \ge a$ implying $b \ge a$, which implies $\sqrt{2} \le 1$. Thus $\sqrt{2} \not\in \Q$. \end{quote}

Since then, several other geometric proofs have been devised for $\sqrt{2}$. Notable examples include the proof by Tom Apostol, based on similar right triangles \cite{Apostol, Bogomolny} and another by Grant Cairns, which involves a different construction using squares \cite{Cairns, Bogomolny}. (The geometric proof of the irrationality of $\sqrt{8}$ can be derived from that of $\sqrt{2}$, as $\sqrt{8}$ is simply twice $\sqrt{2}$.) By extending this method, Steven J. Miller and David Montague \cite{MillerMontague1, MillerMontague2} provided geometric proofs for the irrationality of $\sqrt{5}$ and certain triangular numbers\footnote{The $n$\textsuperscript{th} triangular number, $T_n$, is $1 + 2 + \cdots + n$, and equals $n(n+1)/2$.} such as $\sqrt{3}$, $\sqrt{6}$, and $\sqrt{10}$. Inspired by their work, Ricardo A. Podest\'a gave a geometric proof for $\sqrt{7}$ in \cite{Podesta} using a distinct construction.

We have two goals: to draw attention to these geometric arguments, and to show how far one can push these ideas. The genesis for this paper was a talk given by the second named author at the Math League International Summer Tournament in Summer 2024, where he described these proofs and posed three challenge problems: modify Tennenbaum's construction to prove the irrationality of $\sqrt{7}$, of $\sqrt[3]{2}$, or $\sqrt{6}$ using hexagons (the proofs by Tennenbaum and Miller-Montague use regular $n$-gons for $\sqrt{n}$ with $n \in \{2, 3, 5\}$ but equilateral triangles for $n \in \{6, 10\}$). We\footnote{The first and third named authors were students at the talk, who sent in correct solutions to the third problem and collaborated with each other and the second named author on the other items in the paper. After a draft was written, similar solutions to the third problem were also received by other members of the audience: Cindy Gu (who also provided a nice proof of the irrationality of $\sqrt{7}$, though not by areas, Ningze Song and Roy Sun.} show that a regular $n$-gon argument using hexagons works for $n=6$, and discuss explicitly the limitations in the equilateral triangle method; though we are able to use their ideas to geometrically prove $\sqrt{15} \not\in \Q$, we know their approach must eventually break down for triangular numbers as $T_8 = 36$ is a perfect square! This places us in an interesting position that we know the method works for small $n$ but must eventually break down! We invite the reader to try and find geometric proofs for the irrationality of square-roots of other elements of $\mathcal{I}$, and to let us and the journal know if successful!

%Through a variety of geometric approaches, we compile a comprehensive list of integers—2, 3, 5, 6, 7, 8, 10, and 15—where geometric proofs exist for the irrationality of their square roots.

%%%%%%%%%%%%%%%%%%%%%%%%%%%%%%%%%%%%%%%%%%%%%%%%%%%%%%%%%%%%%%%%%%%%%%%%%%%%%%%%%%%%%%
%%%%%%%%%%%%%%%%%%%%%%%%%%%%%%%%%%%%%%%%%%%%%%%%%%%%%%%%%%%%%%%%%%%%%%%%%%%%%%%%%%%%%%
%%%%%%%%%%%%%%%%%%%%%%%%%%%%%%%%%%%%%%%%%%%%%%%%%%%%%%%%%%%%%%%%%%%%%%%%%%%%%%%%%%%%%%
\section{Geometric proof of the irrationality of $\sqrt{6}$	using hexagons}

\begin{figure}[h]
\begin{center}
\scalebox{.3}{\includegraphics{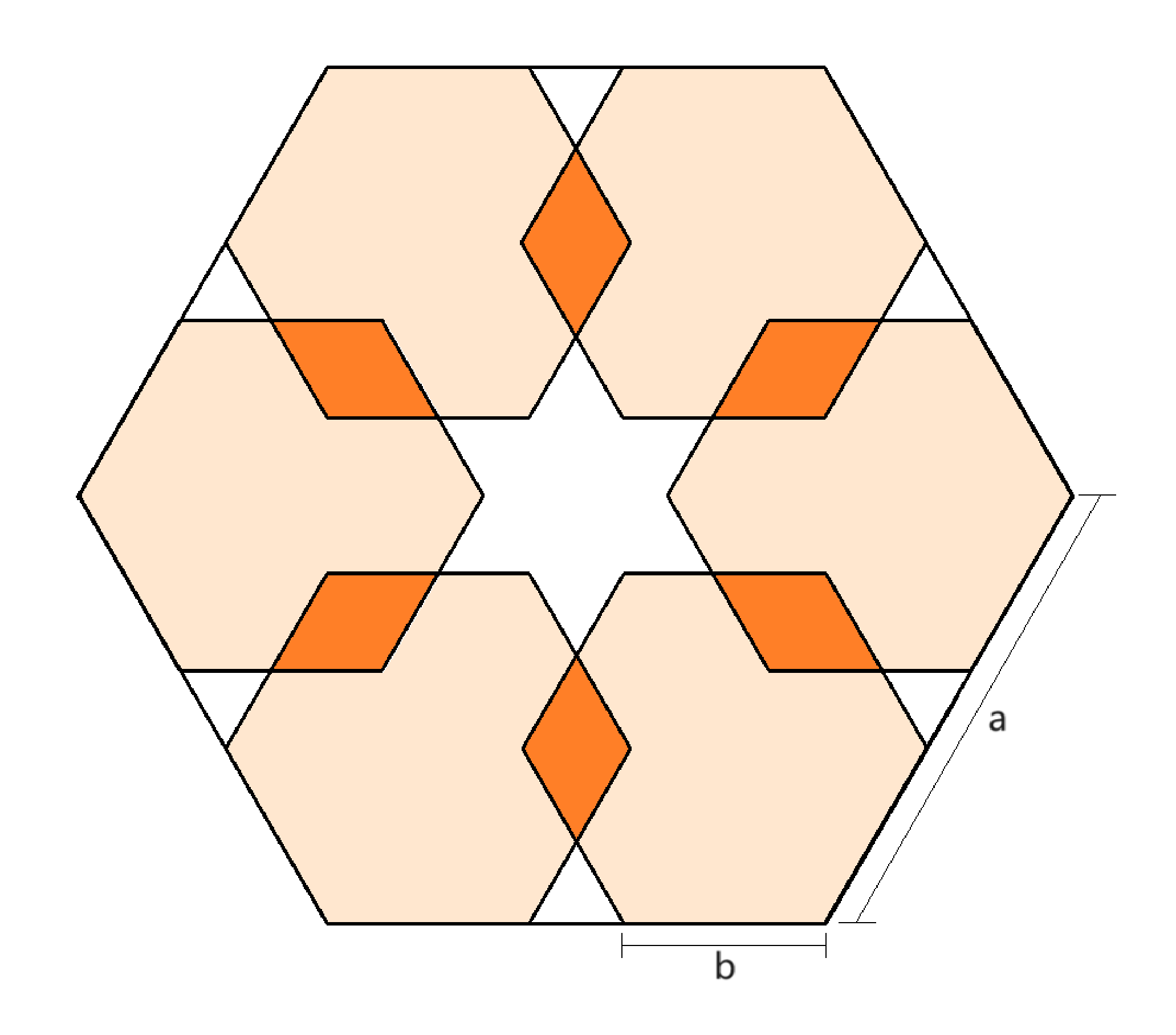}}
\caption{\label{fig:IrrationalitySqrt6hexa} Geometric proof of the irrationality of $\sqrt{6}$}
\end{center}\end{figure}

To prove that \( \sqrt{6} \) is irrational, first suppose \( \sqrt{6} \) is rational and can be expressed as a fraction \( a/b \), where \( a \) and \( b \) are the smallest possible natural numbers satisfying this relation; obviously they are relatively prime or we would have smaller numbers. Then we get \( \left(a/b\right)^2 = 6 \) (or \( a^2 = 6b^2 \)). We interpret \( a^2 = 6b^2 \) geometrically as the area of a large hexagon of length $a$ that equals six times the area of a hexagon of length $b$. We place six regular hexagons of side length \( b \) at the corners of a larger regular hexagon with side length \( a \) (see Figure \ref{fig:IrrationalitySqrt6hexa}).

\begin{figure}[h]
\begin{center}
\scalebox{.2}{\includegraphics{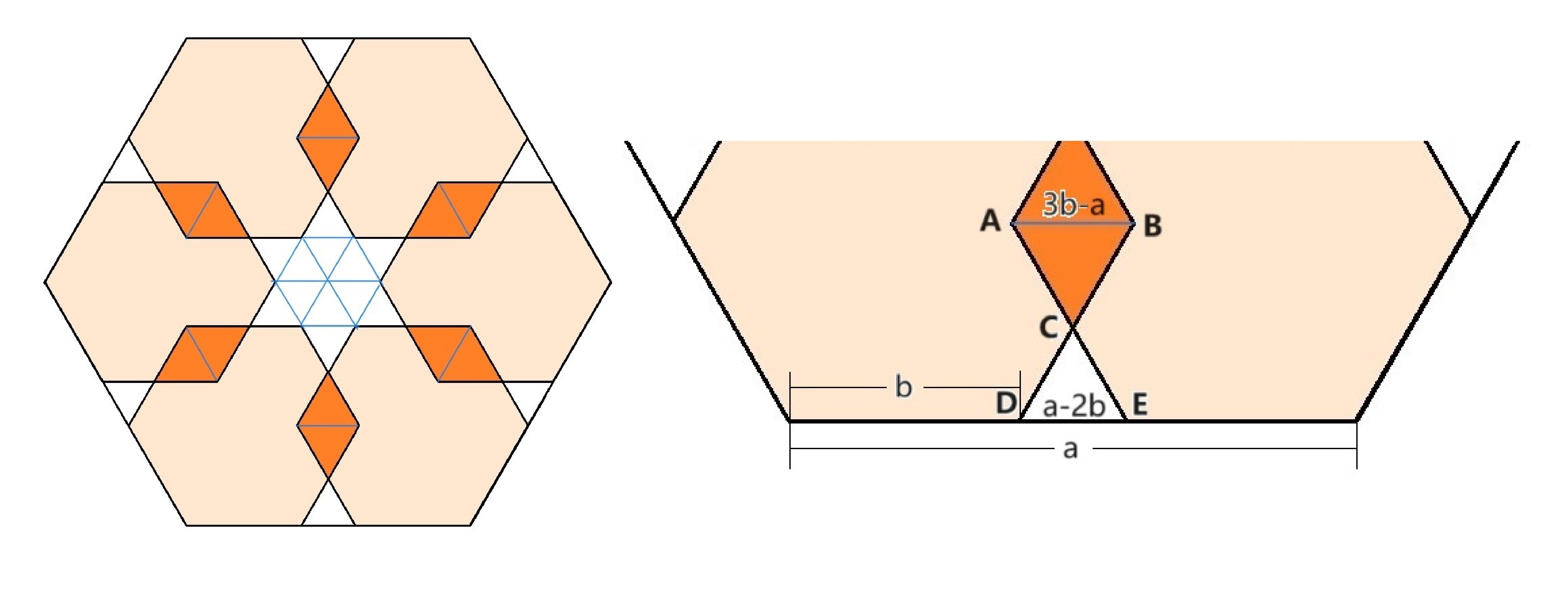}}
\caption{\label{fig:IrrationalitySqrt6hexaLabel} Geometric proof of the irrationality of $\sqrt{6}$: equilateral triangles}
\end{center}\end{figure}

We need to determine the lengths of the sides of the small white and shaded triangles; we clearly label all quantities in Figure \ref{fig:IrrationalitySqrt6hexaLabel} to facilitate the discussion.

First, consider the smallest white triangle \( \Delta CDE \) at the bottom. Its two lower internal angles are both \( \pi/3 \) as they are supplementary angles to the angle of a regular hexagon (\( 2\pi/3 \)), forcing the top angle to be \( \pi/3 \). So, \( \Delta CDE \) is an equilateral triangle. Next, the large white star in the middle is divided into 12 triangles, as shown in \ref{fig:IrrationalitySqrt6hexaLabel}. The 12 smallest white triangles which are each adjacent to an orange diamond are all congruent equilateral triangles by symmetry. This arrangement leaves the hexagon in the center (see the six triangles in the middle of the right image in Figure \ref{fig:IrrationalitySqrt6hexaLabel}) as a regular hexagon composed of six equilateral triangles. Consequently, all 18 of the small white triangles are congruent.

Now divide each overlapping orange diamond into two triangles. In \( \Delta ABC \). The measure of angle ACB is 60 degrees: \( m\angle ACB = \pi/3 \) (for \( \angle ACB \) and \( \angle DCE \) are vertical angles) and \( AC = CB \) (by symmetry). Thus, \( \Delta ABC \) is an equilateral triangle. Additionally, the 12 orange equilateral triangles are congruent by symmetry.

The area of the orange double-counted triangles is equal to the area of the white triangles:
\begin{align*}
12 \times \text{Area (orange equilateral triangle)} &\ = \  18 \times \text{Area (white equilateral triangle)} \\
6 \times \text{Area (orange equilateral triangle)} &\ = \  9 \times \text{Area (white equilateral triangle)} \\
\frac{\text{Area (white equilateral triangle)}}{\text{Area (orange equilateral triangle)}} &\ = \  \frac{6}{9}.
\end{align*}

As the area of an equilateral triangle is proportional to the square of its side length (we don't need this proportionality constant, but it equals $\sqrt{3}/4$), the above relation translates to one about side lengths. Note
\begin{align*}
\frac{\text{Side Length (white equilateral triangle)}}{\text{Side Length (orange equilateral triangle)}} &\ = \  \frac{a-2b}{b-(a-2b)}\ =\ \frac{a-2b}{3b-a}
\end{align*}
and therefore
\begin{align*}
\frac{(a-2b)^2}{(3b-a)^2} &\ = \  \frac{6}{9} \\
\frac{a-2b}{3b-a} &\ = \  \frac{\sqrt{6}}{3} \\
\frac{3a-6b}{3b-a} &\ = \  \sqrt{6}.
\end{align*}

Clearly, \( 3a - 6b \) and \( 3b - a \) are both integers and it's straightforward to show each is positive.\footnote{For example, if $3b - a \le  0$ then $3 \le a/b$ implying $9 \le a^2/b^2 = 6$, a contradiction.} Since \( 2 < \sqrt{6} < 3 \) and $b = a/\sqrt{6}$ we have
\[
3a - 6b\ =\ 3a - \frac{6a}{\sqrt{6}}\ =\ 3a - \sqrt{6}a \ =\ (3 - \sqrt{6})a\ <\ a.
\]
Also using $a = \sqrt{6}b$ we find
\[
3b - a \ = \ 3b - \sqrt{6}b\ =\ (3 - \sqrt{6})b \ <\  b.
\]
We've found a pair of positive integers, \( 3a - 6b \) and \( 3b - a \), that are smaller than \( a \) and \( b \) respectively, with the same ratio of \( \sqrt{6} \). This leads to a contradiction and thus \( \sqrt{6} \) is irrational.

%%%%%%%%%%%%%%%%%%%%%%%%%%%%%%%%%%%%%%%%%%%%%%%%%%%%%%%%%%%%%%%%%%%%%%%%%%%%%%%%%%%%%%
%%%%%%%%%%%%%%%%%%%%%%%%%%%%%%%%%%%%%%%%%%%%%%%%%%%%%%%%%%%%%%%%%%%%%%%%%%%%%%%%%%%%%%
%%%%%%%%%%%%%%%%%%%%%%%%%%%%%%%%%%%%%%%%%%%%%%%%%%%%%%%%%%%%%%%%%%%%%%%%%%%%%%%%%%%%%%
\section{Generalizing Geometric Triangle Irrationality Arguments}

Let's now see how far we can push these geometric arguments. While we know that $\sqrt{n}$ is irrational whenever $n$ is square-free, extending the above method becomes difficult. We invite anyone who can make the regular $n$-gon argument work for additional $n$ to contact us, as we would love to see the answer, but warn you that it is probably hard in general.\footnote{Perhaps if we have a proof that works for $\sqrt{n}$ we can also do $\sqrt{2^k n}$ for any $k$, or maybe if $\sqrt{n_1 n_2}$ if we can do for each $n_i$ and they are relatively prime; we encourage the interested reader to try one of these cases first before hitting the wall, as so many of us have done, with $\sqrt{7}$.} We thus return to the idea introduced by Miller and Montague and look at using multiple equilateral triangles for the triangular numbers, where as always $T_n = n(n+1)/2$.

\begin{figure}[h]
\begin{center}
\scalebox{.1}{\includegraphics{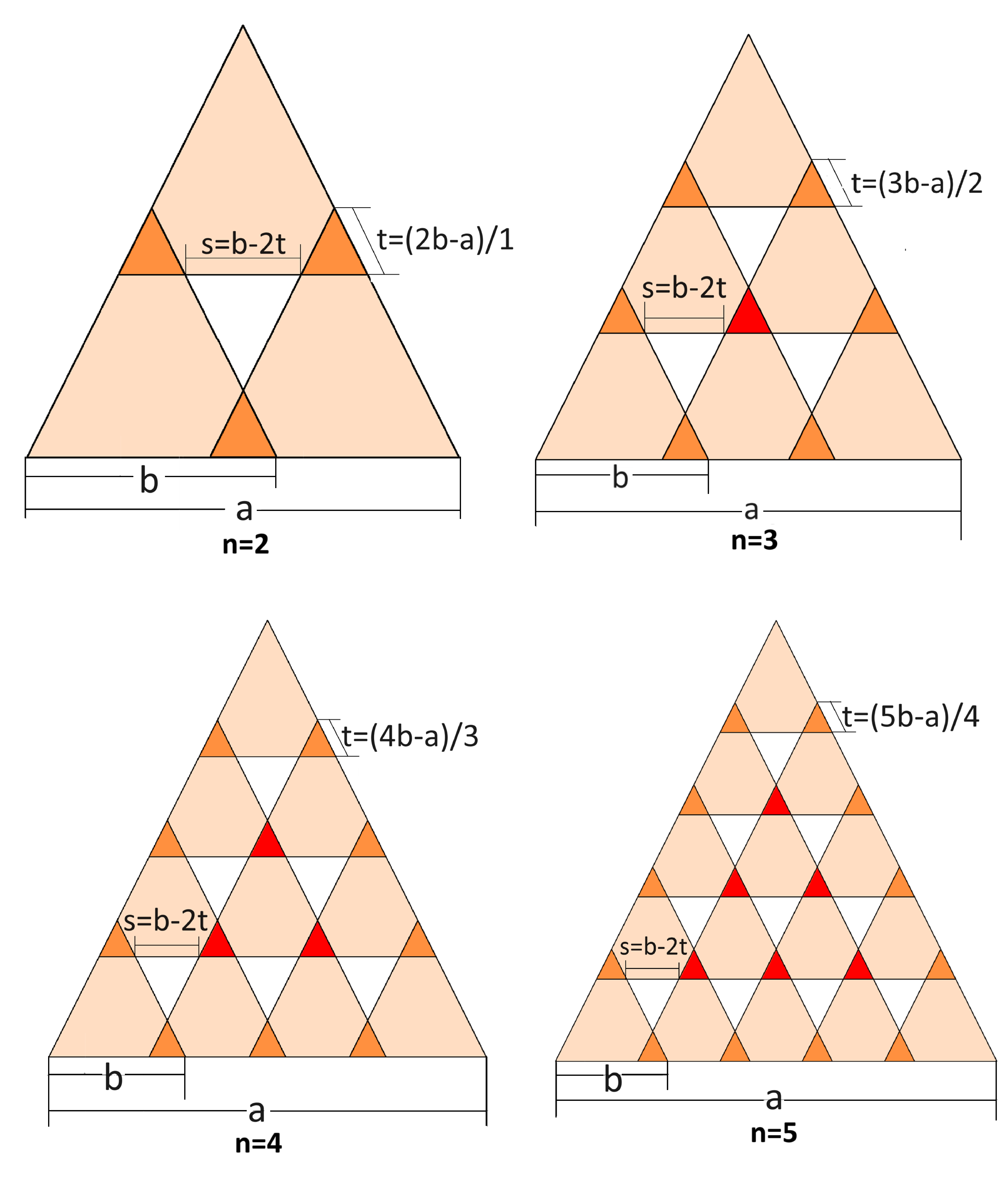}}
\caption{\label{fig:TriangularNumber} The proof of the irrationality of $\sqrt{n(n+1)/2}$ when $n$ is 2, 3, 4, and 5.}
\end{center}\end{figure}

To prove \(\sqrt{n(n+1)/2}\) is irrational, suppose it's rational and thus can be expressed as a fraction $a/b$, where \(a\) and \(b\) are the smallest positive integers that work.\footnote{It's not possible that one fraction has the smallest denominator and another the smallest numerator. If $a_1/b_1 = a_2/b_2$ but $a_1 > a_2$ while $b_1 \le b_2$ then clearly the first fraction exceeds the second.} Thus $$\left(\frac{a}{b}\right)^2 \ =\  \frac{n(n+1)}{2} \ \ \  {\rm or}\ \ \ \ a^2\ =\  \frac{n(n+1)}{2}  b^2.$$ Similar to our earlier arguments, we interpret \( a^2/b^2 = n(n+1)/2 \) geometrically as the ratio of the area of two equilateral triangles. Thus, we arrange \(n\) equally spaced rows of triangles, each with side length \(b\), resulting in a total of \(n(n+1)/2\) triangles inside a larger triangle with side length \(a\). In Figure \ref{fig:TriangularNumber}, four examples of these arrangements are shown for \(n = 2,\ 3,\ 4,\) and \(5\).

\begin{figure}[h]
\begin{center}
\scalebox{.2}{\includegraphics{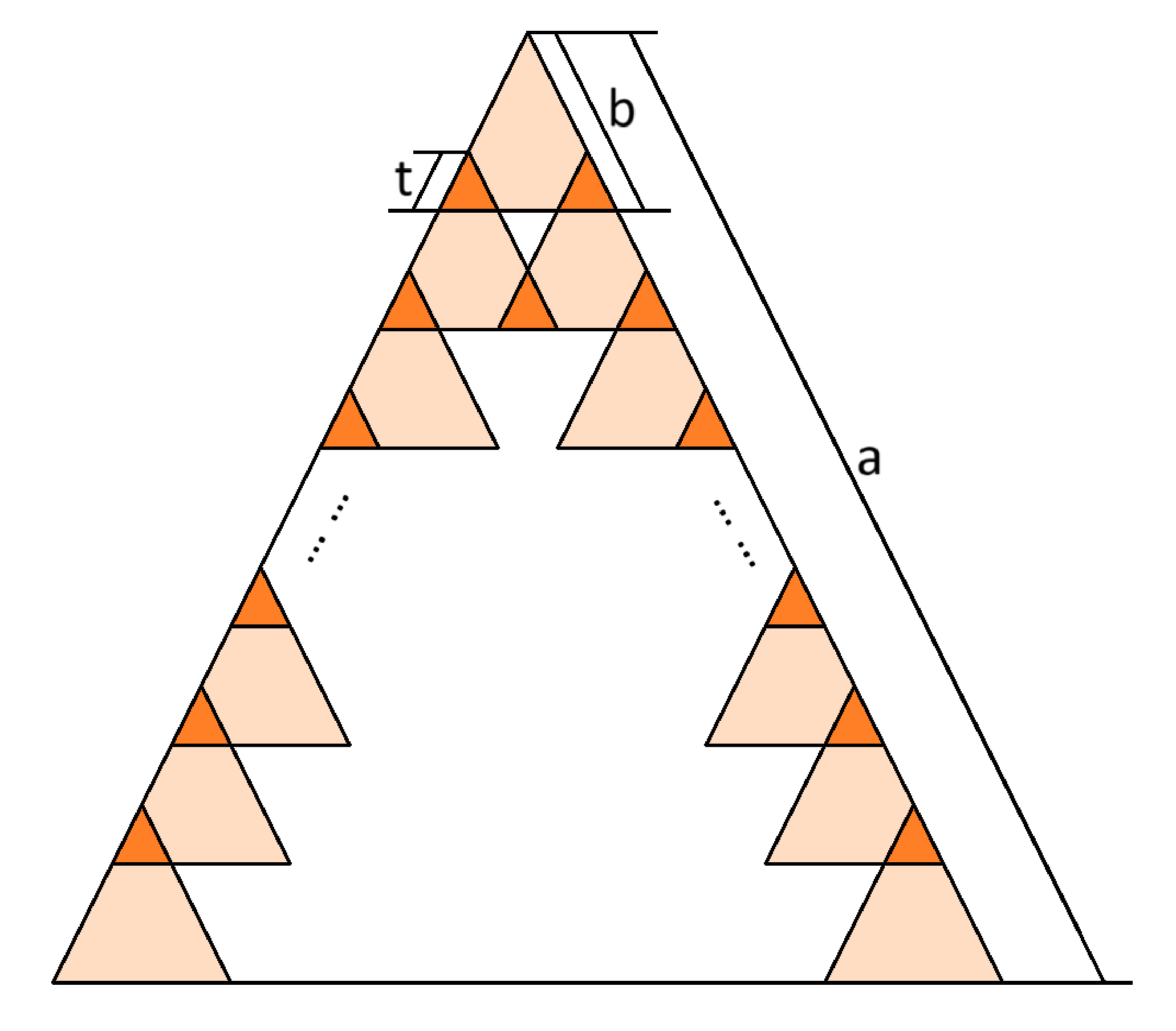}}
\caption{\label{fig:TriangularNumberLegs} Proof of congruence for all doubly-covered and triply-covered triangles: triangles on the outside}
\end{center}\end{figure}

To arrange \(n\) rows of triangles with side length \(b\) with row $k$ having $k$ triangles, we start with the largest triangle with side length \(a\) and place a triangle with side length \(b\) at the top vertex. Next, we place \(n-1\) triangles with side length \(b\) down to the bottom-left vertex and another \(n-1\) to the bottom-right (as shown in Figure \ref{fig:TriangularNumberLegs}). We must be careful in how we do this so that every overlapped triangle is equivalent to every other. By simple continuity arguments we can show this can be done, and give a formula for how far we must move. We start with all $n$ triangles at the top vertex and slide down the left side (and then we mirror this and slide down the right side), moving down a distance $b-t$ and then placing the next triangle, and then moving down another $b-t$ units and placing the next triangle, until we have moved $b-t$ units $n-1$ times and place the final triangle, whose left vertex is now at the left vertex of the triangle of length $a$. Note $t$ is the length of the small equilateral triangles, and each triangle starts $b-t$ units down from the previous. As the final triangle of side length $b$ starts $(n-1)(b-t)$ units down from the vertex and the side of the original triangle is $a$ we must have $$a \ = \ (n-1)(b-t) + b \ \ \ {\rm or}\ \ \ t \ = \ \frac{nb - a}{n-1}. $$

We can conclude that the smallest doubly-covered triangles formed along the legs of the largest equilateral triangle with side length \(a\) are all congruent, as two of their three interior angles correspond to those of an equilateral triangle with side length \(b\), forcing the third to be \(\pi/3\), and all the triangles with side length \(b\) are equidistant. Let the side length of the doubly-covered triangles along the legs of the largest triangle be \(t\). We carefully (perhaps too carefully!) show that all the other doubly and triply counted smaller triangles are equilateral triangles with the same side lengths.

\begin{figure}[h]
\begin{center}
\scalebox{.2}{\includegraphics{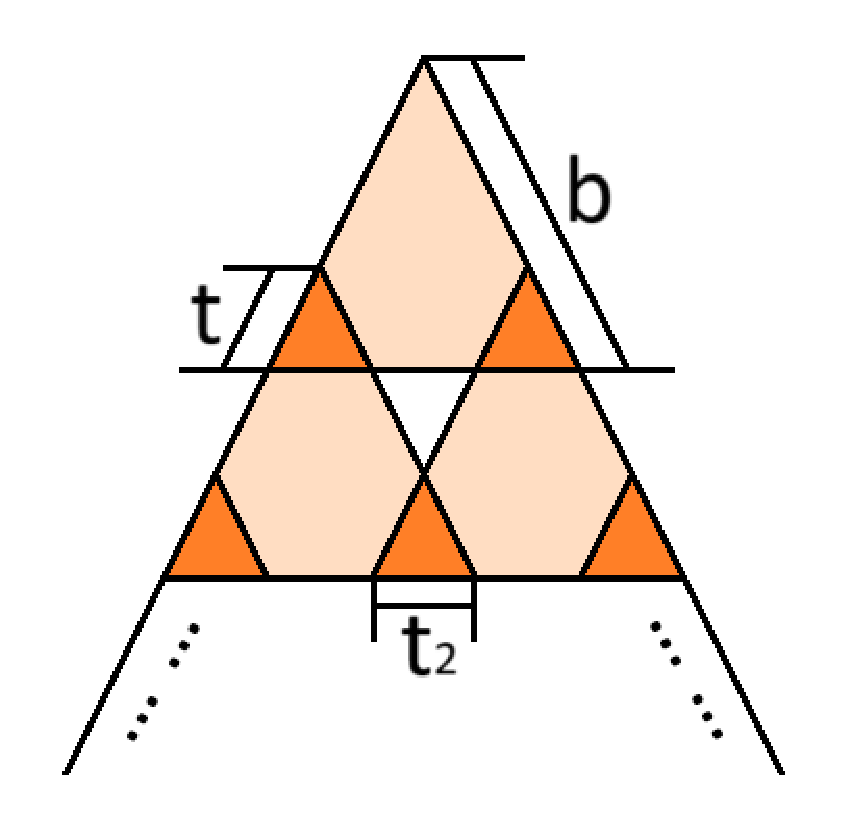}}
\caption{\label{fig:TriangularNumberR2} Proof of congruence for all doubly-covered and triply-covered triangles: first two rows}
\end{center}\end{figure}

First, let us examine the top two rows of triangles with side length \(b\). Let the side length of the overlapping triangle in the middle of the second row at the bottom be \(t_2\) (as shown in Figure \ref{fig:TriangularNumberR2}). The boundary of the left and right sides has a length of \(2b - t\) (two triangles with side length \(b\) subtracting the overlap). The boundary on the bottom has length \(2b - t_2\) (two triangles with side length \(b\) subtracting the overlap). The boundary of the top two rows of triangles with side length \(b\) is an equilateral triangle since each of its interior angles correspond to those of an equilateral triangle with side length \(b\), meaning \(2b - t = 2b - t_2\). Therefore, \(t = t_2\).

\begin{figure}[h]
\begin{center}
\scalebox{.2}{\includegraphics{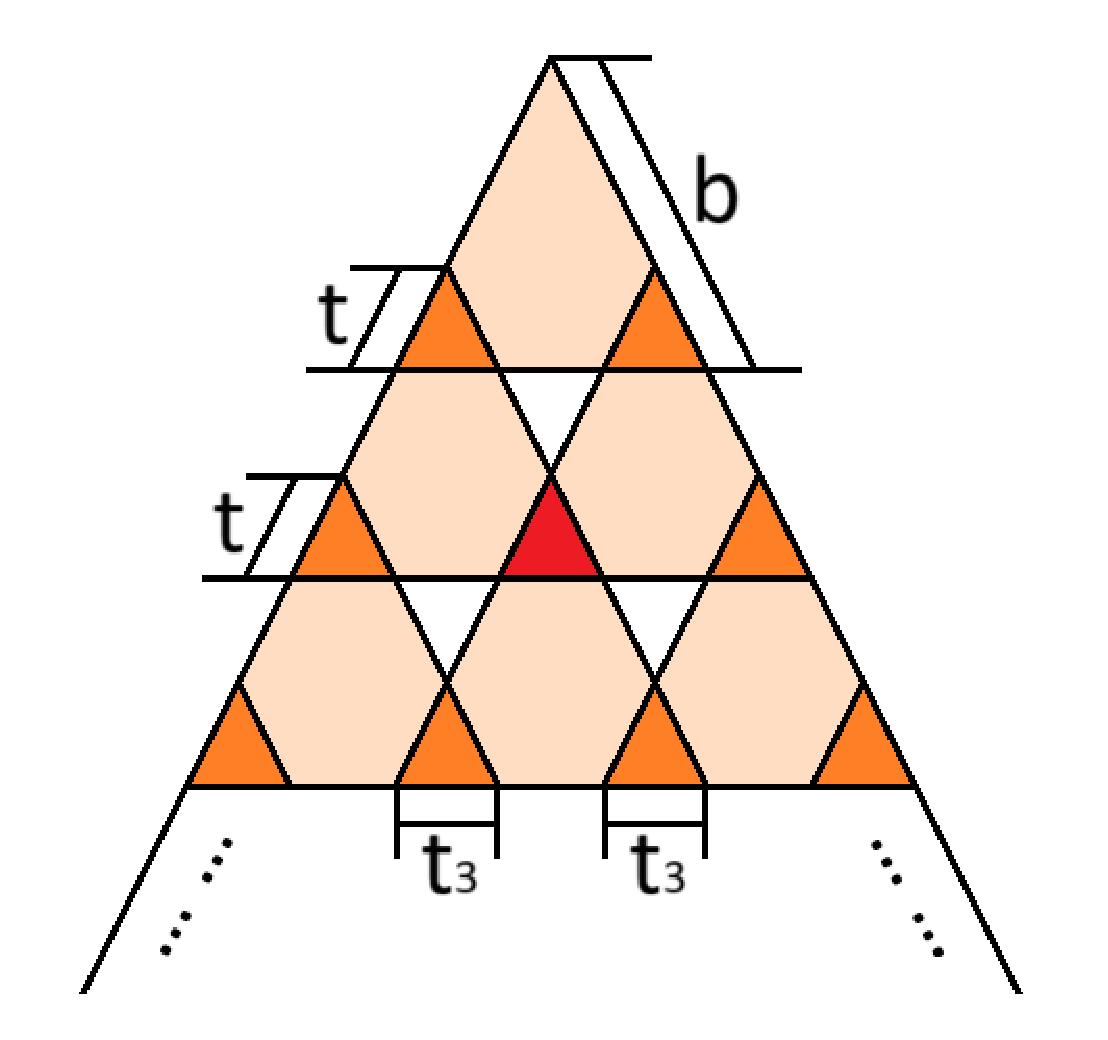}}
\caption{\label{fig:TriangularNumberR3} Proof of congruence for all doubly-covered and triply-covered triangles: first three rows}
\end{center}\end{figure}

Next, we place a triangle with side length \(b\) in the middle of the third row, forming an equilateral triangle as the boundary and adding two more congruent doubly-covered equilateral triangles with side length \(t_3\) on the bottom of the third row (as shown in Figure \ref{fig:TriangularNumberR3}). The boundary of the left and right sides has length \(3b - 2t\) (three triangles with side length \(b\) subtracting the two overlaps). The boundary on the bottom has length \(3b - 2t_3\) (three triangles with side length \(b\) subtracting the two overlaps). Since \(3b - 2t = 3b - 2t_3\), we get \(t = t_3\).

\begin{figure}[h]
\begin{center}
\scalebox{.2}{\includegraphics{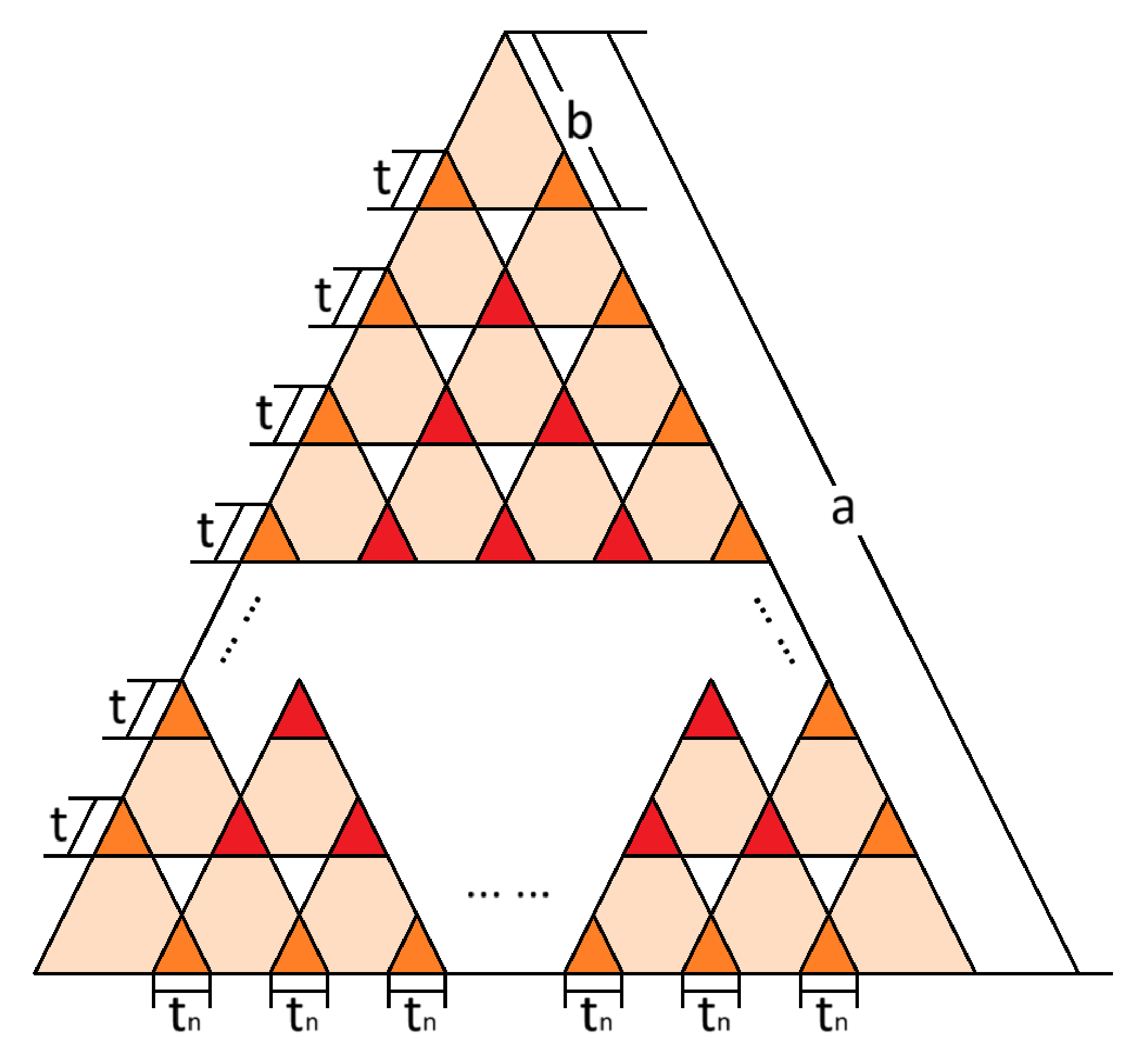}}
\caption{\label{fig:TriangularNumberRn} Proof of congruence for all doubly-covered and triply-covered triangles}
\end{center}\end{figure}

We continue this process until \(n\) rows of triangles of side length \(b\) are arranged within the largest triangle, resulting in multiply-covered equilateral triangles with side length \(t_2, t_3, \ldots, t_n\) in rows 2, 3, \ldots, \(n\) respectively (as shown in Figure \ref{fig:TriangularNumberRn}). The same argument shows each $t_i$ equals $t$, and thus all the triangles are equilateral with the same side length.

%Then, in any row \(k\), \(k-2\) triangles with side length \(b\) are equidistributed in the middle of the \(k\)th row, forming an equilateral triangle as the boundary and adding \(k-1\) more congruent doubly-covered equilateral triangles with side length \(t_k\) on the bottom of the \(k\)th row. The boundary of the left and right sides has length \(kb - (k-1)t\) (\(k\) triangles with side length \(b\) subtracting the \(k-1\) overlaps). The boundary on the bottom has length \(kb - (k-1)t_k\) (\(k\) triangles with side length \(b\) subtracting the \(k-1\) overlaps). Since \(kb - (k-1)t = kb - (k-1)t_k\), we get \(t = t_k\), meaning the multiply-covered triangles on the legs are congruent to the ones on the inside. Therefore, all of the multiply-covered equilateral triangles are congruent and have side length \(t\).

Recall \( t = (nb-a)/(n-1) \). We need to figure out how many doubly and triply counted equilateral triangles of side length $t$ there are, as well as how many triangles are missed (the blank or white triangles) and what their size is. The blank triangles are congruent equilateral triangles since each of their angles is a vertical angle to an interior angle of a multiply covered equilateral triangle and each side of these blank triangles has a length of \( s = b - 2t \) (a triangle with side length \( b \) subtracting two overlapping segments at each end of the side). Then, we have
\[ s \ = \ b - 2t \ = \ b - 2 \cdot \frac{nb - a}{n-1} \ = \ \frac{2a - (n+1)b}{n-1}. \]

The doubly-covered triangles with side length \( t \) are the ones on the outside of the original triangle; there are \( 3(n-1) \), because on each side there are \( n-1 \) doubly-covered triangles. All interior small triangles are triply-covered with side length \( t \). The number of such is \( 1 + 2 + \cdots + (n-2)=(n-2)(n-1)/2 \); note we have $n-1$ rows and $k$ such triangles in row $k$ (when \( n = 2 \) there are no triply-counted triangles). What do these triangles contribute to the doubly and triply counted area? The doubly counted triangles contribute their area and the triply counted ones contribute twice their area (we are looking at the excess coverage, thus one copy of each is needed to get the area of the original triangle of side length $a$, and the area of the extra copies equals the area of the blank triangles). Thus the area they contribute $$\frac{\sqrt{3}}{4}\left[3(n-1) + 2 \frac{(n-2)(n-1)}{2}\right] t^2 \ = \ \frac{\sqrt{3}}{4} (n^2 -1) t^2 \ = \ \frac{\sqrt{3}}{4} (n-1)(n+1) t^2.$$

%Since the area of the total coverings of the smaller triangles with side length \( t \) equals the area of the blank triangles, we need to count the number of each, respectively.

%The grand total of coverings of the smaller triangles with side length \( t \) is the sum of the number of doubly-covered triangles on the outside and two times the number of triply-covered triangles on the inside. Thus, the grand total of coverings of the smaller triangles with side length \( t \) is \[ 3(n-1) + 2[1 + 2 + \cdots + (n-2)] \ = \ (n-1)(n+1). \]

What of the blank triangles? Their area must equal the above. The number of these uncovered triangles with side length \( s \) is \( 1 + 2 + \cdots + (n-1) = (n-1)n/2 \), since, similar to the triply counted triangles, they form a triangular configuration with \( n-1 \) rows. Thus the area of the blank triangles is $\frac{\sqrt{3}}{4} \frac{(n-1)n}{2}s^2$, and setting the two areas equal to each other (and dividing by the scaling factor $\sqrt{3}/4$ yields
\[ (n-1)(n+1) \cdot t^2 \ = \ \frac{(n-1)n}{2} \cdot s^2. \]

We substitute \( t = (nb-a)/(n-1) \) and \( s = [2a-(n+1)b]/(n-1) \) and find
\[ (n-1)(n+1) \cdot \left( \frac{nb - a}{n-1} \right)^2 \ = \ \frac{n(n-1)}{2} \cdot \left( \frac{2a - (n+1)b}{n-1} \right)^2, \]
which by multiplying by $n-1$ (which we can as $n>1$) is equivalent to
\begin{equation}
    (n+1) \cdot (nb - a)^2 \ = \ \frac{n}{2} \cdot [2a - (n+1)b]^2.
    \label{eq:1}
\end{equation}

We now investigate for what $n$ is this solvable. Remember to prove $\sqrt{T_n}$ is irrational we are trying to find a smaller solution to $\sqrt{T_n}$ than $a/b$; this cannot be possible for all $n$ as $T_8 = 36$ has an integer square-root.\\ \

\noindent\textbf{Case 1: $n$ is even.}

Since \( n \) is even, $n/2$ is an integer. We multiply both sides of \eqref{eq:1} by \( n/2 \) to achieve a smaller solution to \( a^2 = n(n+1)/2\cdot b^2 \). Thus we have
\[ \frac{n(n+1) \cdot (nb - a)^2}{2} \ = \ \left( \frac{n}{2} \right)^2 \cdot [2a - (n+1)b]^2, \]
we now divide both sides by \( (nb - a)^2 \) to get
\[ \frac{n^2 \cdot [2a - (n+1)b]^2}{4(nb - a)^2} \ = \ \frac{n(n+1)}{2}, \]
and then taking square-roots yields
\[ \frac{n \cdot [2a - (n+1)b]}{2(nb - a)} \ = \ \sqrt{\frac{n(n+1)}{2}}, \]
and finally since $n$ is even we can move the 2 in the denominator to the numerator:
\[ \frac{n/2 \cdot [2a - (n+1)b]}{nb - a} \ = \ \sqrt{\frac{n(n+1)}{2}}. \]

Since \( n/2 \) is an integer, \( n/2 \cdot [2a - (n+1)b] \) and \( nb - a \) are both integers. Thus to get a contradiction, all we need is the new side \( b' = nb - a \) to be less than \(b\). That is, we need to find all even $n$ with \( nb - a < b \).

As \( a = \sqrt{n(n+1)/2} \cdot b \), such $n$ must satisfy
\[ nb - \sqrt{\frac{n(n+1)}{2}} \cdot b \ = \ \left(n - \sqrt{\frac{n(n+1)}{2}}\right) \cdot b \ < \ 1 \cdot b, \]
which is equivalent to
\[ n - \sqrt{\frac{n(n+1)}{2}} \ < \ 1. \]
Solving the inequality yields  %Solve[n - Sqrt[n (n + 1)/2] == 1, n]
\[ n \ < \ \frac{\sqrt{17} + 5}{2} \ \approx \ 4.56, \]
which is satisfied only for $n=2$ and $n=4$ (recall $n$ is even). Thus, interestingly, even though $T_6 = 21$ has an irrational square-root, this argument breaks down and again we encourage any interested readers who pursue this to reach out to us. \\ \

\noindent\textbf{Case 2: $n$ is odd.}

We multiply both sides of \eqref{eq:1} by \( n+1 \) to achieve a smaller solution to \( a^2 \ = \ n(n+1)/2 \cdot b^2 \), giving us
\[ (n+1)^2 \cdot (nb - a)^2 \ = \ \frac{n(n+1)}{2} \cdot [2a - (n+1)b]^2. \]

Dividing both sides by \( [2a - (n+1)b]^2 \) and taking square-roots yields
\[ \frac{(n+1) \cdot (nb - a)}{2a - (n+1)b} \ = \ \sqrt{\frac{n(n+1)}{2}}. \]

Similar to the case when $n$ was even, we have to group terms appropriately. We pull out a 2 from the denominator and move it to the numerator and have it dividing $n+1$, as $n$ odd implies $n+1$ is even and thus $(n+1)/2$ is an integer. We thus have
\[ \frac{ \frac{n+1}{2} (nb - a)}{a - \frac{n+1}{2} b} \ = \ \sqrt{\frac{n(n+1)}{2}}, \]
with \( \frac{n+1}{2} (nb - a) \) and \( a - \frac{n+1}{2} b \) both integers.

To get a contradiction, we need the new side \( b' = a - \frac{n(n+1)}{2}b \) to be less than \( b \). That is, for what odd $n$ is \( a - \frac{n+1}{2} b < b \)?
We substitute \( a = \sqrt{n(n+1)/2} \cdot b \) and find
\[ a - \frac{(n+1)}{2} b \ = \ \sqrt{\frac{n(n+1)}{2}} b - \frac{(n+1)}{2} b \ = \ \left( \sqrt{\frac{n(n+1)}{2}} - \frac{(n+1)}{2} \right) \cdot b \ < \ 1 \cdot b, \]
which is equivalent to
\[ \sqrt{\frac{n(n+1)}{2}} - \frac{(n+1)}{2} \ < \ 1. \]
Solving the inequality yields  %Solve[Sqrt[(n + 1) n/2] -  (n + 1)/2 == 1, n]
\[ n \ < \ \sqrt{13} + 2 \ \approx \ 5.61, \]
and thus the argument works for odd $n$ of $3$ and $5$ but no greater. \\ \

In conclusion, our geometric method of using equilateral triangles proves the irrationality of the square-root of the triangular numbers\footnote{The ON-LINE ENCYCLOPEDIA OF INTEGER SEQUENCES \cite{OEIS} is a great resource to gain more information on a sequence of triangular numbers that are perfect squares; see \bburl{https://oeis.org/A001108}.} $T_n = n(n+1)/2$ for $n \in \{2, 3, 4, 5\}$ and can be shown to be inapplicable for larger $n$. In particular, $\sqrt{3}, \sqrt{6}, \sqrt{10}$ and $\sqrt{15}$ are all irrational. As remarked many times, it's not unexpected that these triangle games \emph{must} have an obstruction at some point, since there exist infinitely many triangular numbers that are perfect squares; we invite the reader to try these constructions for $n=6$ or $7$ to see what goes wrong.

%%%%%%%%%%%%%%%%%%%%%%%%%%%%%%%%%%%%%%%%%%%%%%%%%%%%%%%%%%%%%%%%%%%%%%%%%%%%%%%%%%%%%%
%%%%%%%%%%%%%%%%%%%%%%%%%%%%%%%%%%%%%%%%%%%%%%%%%%%%%%%%%%%%%%%%%%%%%%%%%%%%%%%%%%%%%%
%%%%%%%%%%%%%%%%%%%%%%%%%%%%%%%%%%%%%%%%%%%%%%%%%%%%%%%%%%%%%%%%%%%%%%%%%%%%%%%%%%%%%%

\ \\

\end{document}